\input amstex
\documentstyle{amsppt}
\magnification 1200
\hsize=14truecm
\hfuzz 1.5pt
\nologo
\NoBlackBoxes

\def\bC{{\Bbb C}}

\def\a{\alpha}

\def \CA {\Cal A}
\def \CU {\Cal U}
\def \la {\lambda}
\def \al {\alpha}
\def \ga {\gamma}
\def \part {\partial}

\topmatter
\title
On algebra generated by Chern-Bott $2$-forms on $\Bbb{SL}_n/{B}$
\endtitle
\author B.~Shapiro$^{\ddag}$, M.~Shapiro$^*$
\endauthor
\headline{\hss{\rm\number\pageno}\hss}
\footline{}
\affil $^\ddag$ Department of Mathematics, University of Stockholm\\
S-10691, Sweden, {\tt shapiro\@matematik.su.se}\\
$^*$ Department of  Mathematics, Royal Institute of Technology\\
Stockholm, S-10044, {\tt mshapiro\@math.kth.se} \\
 \endaffil
\abstract
In this short note we give an explicit presentation of the algebra $\CA_n$
spanned by
 the curvature $2$-forms of the standard Hermitian linear bundles over
$\Bbb{SL}_n/{B}$ as the quotient of the polynomial ring.
 The difference between
 $\CA_n$  and  $H^*(\Bbb{SL}_n/{B})$  reflects the fact that $\Bbb{SL}_n/{B}$
  is not  a symmetric space.  This question was raised by V.~I.~Arnold in
  \cite {Ar}.

\bigskip
{\bf Sur l'alg\`ebre engendr\'ee par les $2$-formes de Chern--Bott sur
$SL_n/B$}
\bigskip
{\it R\'esum\'e:}\ Dans cette note nous donnerons une pr\'esentation
explicite, en tant que quotient de l'anneau polynomial, de l'alg\`ebre $A_n$
engendr\'ee par les $2$-formes de courbure des fibr\'es en droites hermitiens
standard sur $SL_n/B$. La diff\'erence entre $A_n$ et $H^*(SL_n/B)$ refl\`ete
le fait que $SL_n/B$ n'est pas un espace symm\'etrique. Cette question a
\'et\'e
pos\'ee par V.~I.~Arnold dans \cite {Ar}.
\bigskip

\endabstract
\endtopmatter

\document

\subheading {Version fran\c caise abr\'eg\'ee}

Consid\'erons l'espace $SL_n/B$ des drapeaux complets dans $\hbox{\bf C}^n$,
aussi bien que les $n$ fibr\'es quotient $L_i=E_i/E_{i-1}$, o\`u $E_i$
signifie le fibr\'e tautologique standard de dimension $i$ au-dessus de
$SL_N/B$. Quand on fixe une m\'etrique hermitienne dans $\hbox{\bf C}^n$
tous les $L_i$ deviennent des fibr\'es hermitiens. On notera $\omega_i$ la
$2$-forme de courbure de $L_i$. Il est \'evident que chaque $\omega_i$ est
$U_n$-invariante sur $SL_n/B$. Dans cette note nous allons pr\'esenter
l'alg\`ebre $A_n=A(\omega_1,\dots,\omega_n)$, engendr\'ee par
$\omega_1,\dots,\omega_n$, comme un quotient de l'anneau polynomial en $n$
g\'en\'erateurs, et on la comparera ensuite avec $H^*(SL_n/B)$.
\smallskip
Nous \'etudions aussi l'alg\`ebre $A_{k,n}$ engendr\'ee par $k$ parmi les
formes
$\omega_i$. (Sa structure ne d\'epend pas du choix particulier de
$k$-tuple.) Le r\'esultat principal est le suivant.
\smallskip
{\bf Proposition.\ } {\sl L'alg\`ebre $A_{k,n}$ est une alg\`ebre gradu\'ee,
isomorphe \`a $\hbox{\bf C}[x_1,\dots,x_k]/I_{k,n}$, o\`u $I_{k,n}$ est un
id\'eal engendr\'e par $2^k-1$ polyn\^omes de la forme
$$g_{i_1,\dots,i_j}^{(n)}=(x_{i_1}+\dots+x_{i_j})^{j(n-j)+1}.$$
Ici $\{i_1,\dots,i_j\}$ parcourt l'ensemble des sous-ensembles non vides
de $\{1,\dots,n\}$.}
\bigskip
La conjecture suivante est un r\'esultat des exp\'eriments faits sur
l'ordinateur avec le syst\`eme Macaulay.
\smallskip
{\bf Conjecture.}\ {\sl a)\ La dimension de $A_n=A_{n,n}$, en tant qu'espace
vectoriel, est \'egale au nombre total des for\^ets sur $n$ sommets marqu\'es,
et il existe une base monomiale naturelle pour $A_n$, dont les \'el\'ements
sont index\'es par ces for\^ets.

b)\ La dimension totale de $A_{k,n}$ est un polyn\^ome unitaire en $n$ de
degr\'e $k$.}

\heading \S 1. Introduction. \endheading

Let $\Bbb{SL}_n/{B}$ denote the space of complete flags in $\bC^n$. One has the
obvious sequence of tautological bundles $0\subset E_1\subset ... \subset
E_n=E$
 (where $E$ is the trivial $\bC^n$-bundle over $\Bbb{SL}_n/{B}$) and the
corresponding
$n$-tuple of quotient line bundles $L_i=E_i/E_{i-1}$. Fixing some Hermitian
metric
on the original $\bC^n$ one equips every bundle $E_i$ and $L_i$ with the
 induced Hermitian metric. Let $w_i$ denote the curvature form of the above
Hermitian
metric on $L_i$. Each $w_i$ is a $\CU_n$-invariant  $2$-form on
$\Bbb{SL}_n/{B}$ representing
$c_1(L_i)$ in $H^2(\Bbb{SL}_n/{B})$ (the  $\CU_n$-action is provided by the
chosen metric). Denoting $c_1[L_i]$ by $x_i$ one has that
$H^*(\Bbb{SL}_n/{B})=\frac {\bC[x_1,...,x_n]}{(s_1,s_2,...,s_n)}$ where
the denominator is the ideal generated by all elementary symmetric functions
in variables $x_1,...,x_n$, see e.g. \cite {Bo}. Since we have the standard
representative $w_i$ for each $x_i=c_1(L_i)$ it seems natural to study the
algebra
$\CA_n=\CA(w_1,...,w_n)$ generated by all $w_i$'s and compare it to
$H^*(\Bbb{SL}_n/{B})$.
 We will also discuss the algebra $\CA_{k,n}=\CA(w_{i_1},...,w_{i_k})$
generated by any
$k$ of $w_i$'s (the structure of $\CA_{k,n}$ does not depend on a
particular choice of
a $k$-tuple).
The main result of this short note is as follows.

 {\smc 1.1. Proposition.\,} $\CA_{k,n}$ is a graded algebra isomorphic to
 $\frac {\bC[x_1,...,x_k]}{I_{k,n}}$ where
the ideal $I_{k,n}\,$ is generated by the set of $2^k-1$ polynomials of the
form
$$g_{i_1,...,i_j}^{(n)}=(x_{i_1}+...+x_{i_j})^{j(n-j)+1}\eqno(1)$$
where  $\{i_1,...,i_j\}$ runs over the set of all nonempty subsets of the set
 $\{1,...,n\}$.

{\smc 1.2. Example.} The algebra $\CA_3=\CA_{3,3}$ is isomorphic to  $\frac
{\bC[x_1,x_2,x_3]}{I_{3,3}}$
where $I_{3,3}$ is generated by
$x_1^3,x_2^3,x_3^3,(x_1+x_2)^3,(x_1+x_3)^3,(x_2+x_3)^3,
x_1+x_2+x_3$. The Hilbert series of $A_3$ is $(1,2,3,1)$. (For comparison,
the Hilbert
series for $H^*(\Bbb{SL}_3/{B})$ is $(1,2,2,1)$.)

{\smc 1.3. Remark.\,} One has the standard surjective map of algebras
$\pi:\CA_n\to H^*(\Bbb{SL}_n/{B})$.

The following conjecture is the result of calculation of the Hilbert series
of $\CA_n$ for $n\le 7$, see \cite {SS}.

{\smc 1.4. Conjecture.} 1) The total dimension of $\CA_n$ (as a vector
space) equals the
number of forests on $n$ labeled vertices and there exists a natural
monomial basis
for $\CA_n$ whose monomials are enumerated by the above forests.
\newline
2) The total dimension
of $\CA_{k,n}$ is a monic polynomial in $n$ of degree $k$.

Besides purely aesthetic reasons the study of $\CA_{n}$ is motivated by the
fact that many known results  for $H^*(\Bbb{SL}_n/{B})$ (such as e.g. the
existence
of a good monomial basis, $S_n$-module structure, the action of divided
difference
operators etc) have natural counterparts for $\CA_{n}$. One might hope to
develop
 a theory of characteristic classes with values
in $\CA_n$ for $\CU_n$-invariant Hermitian bundles over $\Bbb{SL}_n/{B}$.
Finally,
the study of $\CA_n$ seems to be important for the understanding of
structure for the
algebra $\frak A (\Bbb{SL}_n/{B})$ of all $\CU_n$-invariant forms on
$\Bbb{SL}_n/{B}$ which recently appeared in the arithmetic
 intersection theory for flag varieties, comp. \cite {Ta1}, \cite {Ta2}.

The authors are  grateful to V.~ Arnold for the formulation of the problem
and to R.~Fr\"oberg for his help with Macaulay which was crucial in
formulation of 1.4.
 Sincere thanks are due to R.Stanley who kindly explained to the first
author some details about counting trees and forests and calculated the
generating function
 for the number of forests on labeled vertices. Finally, the authors want
to thank
H.Tamvakis for the explanation of his papers \cite {Ta1} and \cite {Ta2} as
well as
for some relevant references.

\heading \S 2. Proofs. \endheading

\subheading { On algebra $\frak A(\Bbb{SL}_n/{B})$ of  $\CU_n$-invariant
forms on $\Bbb{SL}_n/{B}$}

One knows that $\Bbb{SL}_n/{B}$  has another presentation as a
homogeneous space, namely $\CU_n/T^n$ where $T^n\subset \CU_n$ is the usual
torus of
diagonal matrices. Let us recall an old result from the general
theory of homogeneous spaces, see e.g. \cite {DNF}.

{\smc 2.1. Proposition.} The ring of $G$-invariant differential forms on a
homogeneous
space $G/H$ is isomorphic to the exterior algebra  $\Lambda_{inv}((\frak
G/\frak H)^*)$,
i.e. to the algebra of skewsymmetric polylinear functions on $\frak G$
which a) vanish
on $\frak H$ and b) are invariant under the action of internal
automorphisms by the
elements of $H$. (Here $\frak G$ and $\frak H$ are the Lie algebras of $G$
and $H$ resp.)

Take the $\binom n 2$-dimensional vector space $V_n$ of all skew-hermitian
matrices of the
form
$$V_n=\pmatrix 0 & a_{1,2} &\dots&\dots  &\dots& a_{1,n}\\
     -\bar a_{1,2}& 0   &a_{2,3}&\dots & \dots& a_{2,n}\\
     -\bar a_{1,3}&-\bar a_{2,3}&0&a_{3,4}&\dots& a_{3,n}\\
      \vdots     &\vdots       &\vdots&\ddots &\vdots&\vdots\\
     -\bar a_{1,n}&-\bar a_{2,n}&\dots&\dots&\dots&0
\endpmatrix . \eqno(2)$$

Its matrix entries form the linear space $(\frak U_n/\frak T^n)^*$
where $\frak U_n$ (resp. $\frak T^n$) is the Lie algebra of $\CU_n$ (resp.
of $T^n$).
Let us denote by
$(e^{i\la_1},e^{i\la_2},...,e^{i\la_n})$ the diagonal entries of elements in
 $T^n$  acting   by conjugation on $V_n$. Under this action each entry
$a_{i,j}$ above the main diagonal is multiplied
by $e^{i(\la_i-\la_j)}$ and each entry $-\bar a_{i,j}$ below the main
diagonal is multiplied
by $e^{i(\la_j-\la_i)}$. Introducing the principal weights
$\al_i=\la_i-\la_{i+1},
i=1,...,n-1$ we get that $a_{i,j}$ is multiplied by
$e^{i(\al_i+\al_{i+1}+...+\al_{j-1}) }$.
The expression $\al_i+\al_{i+1}+...+\al_{j-1}$ is called {\it the
multiweight} of the
entry $a_{i,j}$. (Under this convention the entry $-\bar a_{i,j}$ has the
opposite multiweight
$-\al_i-\al_{i+1}-...-\al_{j-1}$.)  The multiweight of an exterior monomial
$\tilde a_{i_1,j_1}\wedge \tilde a_{i_2,j_2}\dots \wedge \tilde
a_{i_r,j_r}$ where each
$\tilde a_{i_l,j_l}$ is either $a_{i_l,j_l}$ or $-\bar a_{i_l,j_l}$ equals
to the
sum of multiweights of its factors.

{\smc 2.2. Corollary.} The ring $\frak A(\Bbb{SL}_n/{B})$ is the linear
span of all exterior
monomials $\tilde a_{i_1,j_1}\wedge \tilde a_{i_2,j_2}\dots \tilde
a_{i_r,j_r}$ having
vanishing multiweight. (In particular, $\frak A(\Bbb{SL}_n/{B})$ has no
degree 1 elements.)

{\smc 2.3. Example.} The Hilbert series for $\frak A(\Bbb{SL}_3/{B})$ is
$(1,0,3,2,3,0,1)$. (We
assume that $\frak A(\Bbb{SL}_n/{B})$ contains constants.) Its degree 2
component is spanned
by $a_{1,2}\wedge\bar a_{1,2}, a_{1,3}\wedge\bar a_{1,3},a_{2,3}\wedge\bar
a_{2,3}$;
degree 3 component is spanned by $a_{1,2}\wedge\a_{2,3}\wedge \bar a_{1,3},
\bar a_{1,2}\wedge\bar \a_{2,3}\wedge  a_{1,3}$; degree 4 component is
spanned by
$a_{1,2}\wedge\bar a_{1,2}\wedge a_{1,3}\wedge\bar a_{1,3},
a_{1,2}\wedge\bar a_{1,2}\wedge a_{2,3}\wedge\bar a_{2,3},
a_{1,3}\wedge\bar a_{1,3}\wedge a_{2,3}\wedge\bar a_{2,3}$ and, finally,
its degree 6
component is spanned by $a_{1,2}\wedge\bar a_{1,2}\wedge a_{1,3}\wedge\bar
a_{1,3}
\wedge a_{2,3}\wedge\bar a_{2,3}$. The Hilbert series for $\frak
A(\Bbb{SL}_4/{B})$
is $(1,0,6,4,18,12,26,12,18,4,6,0,1)$.

Recall that {\it an Eulerian digraph} is a digraph  such that the numbers
of coming and leaving
edges at each vertex are equal. (We do not allow loops.)

The following proposition is a relatively simple reformulation of 2.2.

{\smc 2.4. Proposition,} see \cite {SS}.  The dimension of $\frak
A(\Bbb{SL}_n/{B})$ equals
$2^{\binom {n}{2}} +2\sharp_{eul}(n)$ where $\sharp_{eul}(n)$ is the number
of Eulerian
digraphs on $n$ labeled vertices.

\subheading {On curvature forms}

Using the above desription of $\CU_n$-invariant forms on $\Bbb{SL}_n/{B}$ one
can be easily derive the following proposition  from the results of \cite
{GS}, see also
\cite {Ta2}.

{\smc 2.5. Proposition.} The curvature $2$-form $w_i,\,i=1,...,n$ of the
tautologic
line bundle $L_i=E_i/E_{i-1}$ over $\Bbb{SL}_n/{B}$ equals to the sum of
all entries
in the $i$th row of the following matrix of $2$-forms

$$\pmatrix 0&\ga_{1,2}&\ga_{1,3}&...&\ga_{1,n}\\
         -\ga_{1,2}&0&\ga_{2,3}&...&\ga_{2,n}\\
         -\ga_{1,3}&-\ga_{2,3}&0&...&\ga_{3,n}\\
         \vdots&\vdots&\vdots&\ddots&\vdots\\
        -\ga_{1,n}&-\ga_{2,n}&-...&-\ga_{n-1,n}&0\endpmatrix \eqno(3)$$
where $\ga_{i,j}=\al_{i,j}\wedge \bar\al_{i,j}$.

Now we are ready to start proving the main proposition 1.1.
Consider the algebra $\CA(w_{i_1},...,w_{i_k})$ generated by
$w_{i_1},...,w_{i_k}$. Since all $w_j$ are $\CU_n$-invariant one has
that $\CA(w_{i_1},...,w_{i_k})$ is a subalgebra of the exterior algebra
$\Lambda(\Bbb C^{\binom n 2})$.

{\smc 2.6. Lemma.} There exists a natural $S_n$-action on the set of all
permutations
$(w_{i_1},...,w_{i_n})$ of $w_1,...,w_n$.

{\smc Proof.\,} The elementary transposition $\tau_i,\, i=1,...,n-1$ acts
on the matrix $(3)$ of $2$-forms by the simultaneous interchanging of
1) the $i$th row with the $(i+1)$st row, 2) the $i$th column with the
$(i+1)$st column and 3)changing sign of $\ga_{i,i+1}$. One can easily check
that
this determines the required $S_n$-action. \qed

{\smc 2.7. Corollary.} All algebras $\CA(w_{i_1},...,w_{i_k})$ are
isomorphic to
each other and, in particular, to $\CA(w_{1},...,w_{k})$.

 We denote this class of isomorphic algebras by $\CA_{k,n}$.

{\smc 2.8. Lemma.\,} $\CA_{n-1,n}$ is isomorphic to $\CA_n=\CA_{n,n}$.

{\smc Proof.\,} Indeed, one has $w_1+w_2+...+w_n=0$. \qed

{\smc Definition.} Let us call by {\it the vanishing ideal}
$I_{k,n}$ of $\CA_{k,n}$ the set of all polynomials $p\in \Bbb C[x_1,...,x_k]$
 which vanish if we substitute the variables $x_1,...,x_k$ by the curvature
forms $w_1,....,w_k$ resp.

{\smc 2.9. Lemma.\,} The vanishing ideal $I_{k,n}$ consists of all $p\in
\Bbb C[x_1,...,x_k]$
such that $2^n$ derivatives  $p, \frac {\part p}{\part x_1},
\frac {\part p}{\part x_2},...,\frac {\part p}{\part x_k},
\frac {\part^2 p}{\part x_1 \part x_2},..., \frac {\part^2 p}{\part x_{k-1}
 \part x_k}, ..., \frac {\part^k p}{\part x_1 \part x_2 ...\part x_k}$
belong to
$I_{k,n-1}$. (Derivatives are taken with respect to all subsets
$\{i_1,...,i_l\}\subset \{1,...,n\}$ of indices (no repetitions of indices
are allowed) including the empty set.)

{\smc Proof.\quad} Notice that  $\CA_{k,n}$ has a natural module structure over
$\CA_{k,n-1}$. Indeed, one has $w_i^{(n)}=w_i^{(n-1)}+\ga_{i,n}$ where the
upper
index shows the dimension of the initial space, see $(3)$. Therefore for
any polynomial
$p\in \Bbb C[x_1,...,x_k]$ one has after substitution of the curvature forms

$ p(w_1^{(n)},...,w_k^{(n)})=
p(w_1^{(n-1)}+\ga_{1,n},...,w_k^{(n-1)}+\ga_{k,n})=
 p(w_1^{(n-1)},...,w_k^{(n-1)})+
 p_{x_1}(w_1^{(n-1)},...,w_k^{(n-1)})\ga_{1,n}+...+p_{x_k}(w_1^{(n-1)},...,w
_k^{(n-1)})
\ga_{k,n}+...+
\quad\quad\quad\quad$ $
p_{x_1,x_2,...,x_k}(w_1^{(n-1)},...,w_k^{(n-1)})\ga_{1,n}\ga_{2,n}...\ga_{k,
n}.
\quad\quad\quad\quad\quad\quad\quad\quad\quad\quad\quad\quad(4) $

Since $\ga^2_{i,n}=0$ the occuring monomials in the r.h.s. are square-free and
the coefficient at the product $\ga_{i_1,n}...\ga_{i_l,n}$ equals
 $\frac {\part^l p}{\part x_{i_1}...\part x_{i_l}}$. Finally, the condition
 $p\in I_{k,n}$, i.e. $p(w_1^{(n)},...,w_k^{(n)})= 0$ is equivalent to
vanishing of all
 polynomial coefficients $p_{x_{i_1},...,x_{i_l}}(w_1^{(n-1)},...,w_k^{(n-1)})$
in the r.h.s of $(4)$. By definition this means  that
$p_{x_{i_1},...,x_{i_l}}(x_1,...,x_k)\in
I_{k,n-1}$.\qed

Denote $\Cal D_{i_1,...,i_l}=\frac{\part^l}{\part x_{i_1} ... \part
x_{i_l}}$ and
let $V_{i_1,...,i_l,r,m,k}$ be the space of all homogeneous polynomials of
degree $m+r$ in $k$
variables of the form $p=(x_{i_1}+...x_{i_l})^rf$ where $deg(f)=m$.

{\smc 2.10. Lemma.\,} The linear map $\Cal D_{i_1,...,i_l}:
V_{i_1,...,i_l,r+l,m+k-l,k}\to
V_{i_1,...,i_l,r,m,k}$ is a surjection.

{\smc Proof.\quad} Simple linear algebra. \qed

\subheading {Proof of the main proposition 1.1}\medskip

 We  use a double induction on $k\le n$, i.e. for
a given $k$ we apply induction on $n$ and then cover the first case
$\CA_{k+1,k+1}$ for $k+1$
using lemma 2.8.

Base of induction. $\CA_{1,n}=\frac {\Bbb C [x_1]}{(x_1^n)}$. Indeed,
$w_1=\ga_{1,n}+\ga_{1,3}+
...\ga_{1,n}$ where $\ga_{1,i}$ are independent Grassmann variables, i.e.
$\ga_{1,i}^2=0$.
Statement follows.

Step of induction. Assume that 1.1. is proven for all pairs $(k^\prime,
n^\prime)$
where $k^\prime<k$ and for all $(k,n^\prime)$ where $n^\prime<n$. Let us
show that it holds for
the pair $(k,n)$. Notice that all polynomials $(1)$ are contained in
$I_{k,n}$. This can be
either checked directly or by induction using 2.9. Let us temporarily denote by
$\tilde I_{k,n}$ the ideal generated by all polynomials in $(1)$. We have
to show that
$I_{k,n}=\tilde I_{k,n}$. By the above remark  $\tilde I_{k,n}\subset I_{k,n}$.
 Take any $p\in I_{k,n}$. Using lemma 2.10 we will step by step add to $p$
certain
polynomials from $\tilde I_{k,n}$ in such a way that the resulting
polynomial $p_{fin}$
will have of derivatives from lemma 2.9 vanishing.
Consider $\Cal D_{1,...,k}(p)$. By 2.9. it belongs to $I_{k,n-1}$. By the
inductive assumption
$I_{k,n-1}$ is generated by $g_{i_1,...,i_l}^{(n-1)}$, see $(1)$. Therefore
one has
$\Cal
D_{i_1,...,i_l}(p)=\sum_{i_1,...,i_l}(x_{i_1}+...+x_{i_l})^{j(n-j-1)+1}h_{i_
1,...,i_l}$.
By lemma 2.10 for each
$\phi=(x_{i_1}+...+x_{i_l})^{j(n-j-1)+1}h_{i_1,...,i_l}$ there exists
(but nonunique) $\psi=(x_{i_1}+...+x_{i_l})^{j(n-j)+1}H_{i_1,...,i_l}$ such
that
$\Cal D_{i_1,...,i_l}(\psi)=\phi$. Notice that $\psi\in \tilde I_{k,n}$.
Therefore subtracting
from $p$ an appropriate polynomial belonging to $\tilde I_{k,n}$  we get
$\tilde p\in I_{k,n}$ such
that $\Cal D_{1,...,k}(\tilde p)=0$. Consider now any derivative
$\Cal D_{1,...,\hat i,...,k}=\frac{\part^k}{\part x_{1} ...
\widehat{\part x_{i}}...\part x_k}$  applied to $\tilde p$.
One has that $\Cal D_{1,...,\hat i,...,k}(\tilde p)$ does not
depend on $x_i$ since $\Cal D_{1,...,k}(\tilde p)=0$. Therefore using the
same argument as
above we can subtract from $\tilde p$ some polynomial belonging to $\tilde
I_{k,n}$ which does not
depend on $x_i$ and such that $\Cal D_{1,...,\hat i,...,k}$ applied to the
resulting
difference vanish. Notice that if we apply this procedure for each $\Cal
D_{1,...,\hat i,...,k}$
consecutively to the polynomial obtained on  the previous step we will not
change vanishing
of all previous $\Cal D_{1,...,\hat j,...,k},\,j<i$ apllied to the
polynomial obtained on the
current step since we add a function which does not depend on $x_i$. Having
obtained vanishing
of all  $\Cal D_{1,...,\hat i,...,k},\,i\le k$ applied to our polynomial we
can proceed with all
 $\Cal D_{1,...,\hat i,...,\hat j,...,k}$ (ordering them e.g.
lexicographically). Our
assumptions imply that $\Cal D_{1,...,\hat i,...,\hat j,...,k}$ applied to
our polynomial
does not depend on $x_i$ and $x_j$. Therefore we can subtract a polynomial
from $\tilde
I_{k,n}$ which does not depend on $x_i$ and $x_j$ either and such that
$\Cal D_{1,...,\hat i,...,\hat j,...,k}$
applied to the resulting difference vanish. Using this procedure
consecutively we do not
spoil vanishing of the previous derivatives since the polynomials we
subtract depend on different
groups of variables. Continuing in the same manner we get some polynomial
$p_{fin}$  all
derivatives of which mentioned in 2.9 vanish.  This means that $p_{fin}$ is
a constant.
But since constants different from zero do not lie in $I_{k,n}$ we have that
$p_{fin}$ should also vanish. The statement follows.
\qed

\heading \S 3. Some related problems. \endheading

{\smc Problem 1.\,} Give a presentation of the obvious analog
of $\CA_n$ for any $\Bbb{SL}_n/{P}$ generated by the standard $\CU_n$-invariant
forms representing  the Chern classes of the tautological (quotient) bundles.

{\smc Remark.\,} Notice that in the case of Grassmanian $G_{k,n}$ the analogous
algebra coincides with $H^*(G_{k,n})$ since $G_{k,n}$ is a symmetric
space and therefore has no nontrivial left-invariant forms homologous to
$0$, see e.g. \cite {St}.

{\smc Problem 2.\,} Prove conjecture 1.4 and determine the $S_n$-module
structure for $\CA_n$ and its
$\Bbb{SL}_n/{P}$-analogs.

{\smc Problem 3.\,} Calculate the Poincare series for the algebra $\frak
A(\Bbb{SL}_n/{P})$
of $\CU_n$-invariant forms on $\Bbb{SL}_n/{P}$.


\Refs
\widestnumber\key{SSV1}

\ref \key Ar \by V.~I.~Arnold
\paper Remarks on eigenvalues and eigenvectors of Hermitian matrices,
Berry phase, adiabatic connections and quantum Hall effect
\jour Selecta Math.
\yr 1995
\pages 1--19
\vol 1
\issue 1
\endref

\ref \key Bo \by A.~Borel
\paper La cohomologie mod 2 des certains espaces homog\`enes
\jour Comment. Math. Halv
\yr 1953
\pages 165--197
\vol 27
\endref

\ref \key DNF \by B.~Dubrovin, S.~Novikov and A.~Fomenko
\book Modern geometry
\vol 2
\publ Moscow, Nauka
\yr 1984
\endref

\ref \key GS \by P.~Griffiths and W.~Schmid
\paper Locally homogeneous complex manifolds
\jour Acta Math
\vol 123
\yr 1969
\pages 253--302
\endref


\ref \key MSh \by M.~Shapiro and A.~ D.~ Vainshtein
\paper Stratification of Hermitian matrices,  the Alexander mapping,
and the bundle of eigenspaces
\jour Compt. Rendus
\yr 1995
\vol 321
\issue 12
\pages 1599--1604
\endref

\ref \key SS \by B.~Shapiro and M.~Shapiro
\paper Algebra of Chern-Bott forms on $\Bbb{SL}_n/{B}$  and
labeled forests
\finalinfo in preparation
\endref


\ref \key St \by W.~Stoll
\book Invariant forms on Grassmann manifolds
\publ Annals of Math. Studies, Princeton Univ. Press
\yr 1977
\pages 113
\vol 89
\endref

\ref \key Ta1 \by H.~Tamvakis
\paper Bott-Chern forms and arithmetic intersections
\finalinfo alg-geom/9611005, nov. 96
\endref

\ref \key Ta2 \by H.~Tamvakis
\paper Arithmetic intersection theory on flag varieties
\finalinfo alg-geom/9611006, nov. 96
\endref


\endRefs
\enddocument